\documentclass[10pt]{article}
\usepackage[cp1251]{inputenc}
\usepackage[russian]{babel}
\usepackage{amsmath}
\usepackage{amssymb}
\usepackage{amsthm}
\usepackage{latexsym}

\theoremstyle{plain}

\begin{document}
\begin{center}
{\bf Inversion formula for infinitely divisible distributions}
\end{center}
\begin{center}
{\bf E.V. Burnaev}
\end{center}

The aim of this note is to prove the inversion formula, which can
be used to compute the Levi measure of an infinitely divisible
distribution from its characteristic function. Obtained formula is
similar to the well-known inversion formula [2], which is used to
compute the distribution function of a random variable from
the corresponding characteristic function.\\

{\sc Theorem.} {\em Let $\phi(z)=\mathrm{e}^{\psi(z)}$ be the
characteristic function of an infinitely divisible distribution
with the cumulant function $\psi(z)=i\gamma
z-1/2\cdot\sigma^2z^2+\int_{\mathbb{R}}\left(\mathrm{e}^{izx}-1-izxh(x)\right)U(dx)$,
where $h(x)$ is a truncation function and $U(dx)$ is the Levy
measure. Let
$\hat{\rho}(z)=\int_{-1}^1\log\frac{\phi(z)}{\phi(z+\lambda)}d\lambda-1/3\cdot\sigma^2$
 and $B=[a,b]$, where $-\infty<a<b<\infty$ and $U(\{a\}\cup\{b\})=0$. We then have
\begin{equation}
\int_B\left(\frac{x-\sin
x}{x}\right)U(dx)=\frac{1}{4\pi}\cdot\mathrm{v.p.}\int_{\mathbb{R}}
\hat{\rho}(z) \left(\int_B\mathrm{e}^{-izx}dx\right)dz. \label{1}
\end{equation}
Moreover, if $\int_{|x|>1}x^2U(dx)<\infty$, then
\begin{equation}
\int_Bx^2U(dx)=-\frac{1}{2\pi}\cdot\mathrm{v.p.}\int_{\mathbb{R}}
\left(\psi''(z)+\sigma^2\right)\left(\int_{B}\mathrm{e}^{-izx}dx\right)dz.
\label{2}
\end{equation}
If $\int_{\mathbb{R}}\left|\hat{\rho}(z)\right|dz<\infty$, then
the Levy measure is absolutely continuous with respect to the
Lebesgue measure, has a bounded continuous density $u(x)$, and
\begin{equation} u(x)=\frac{x}{x-\sin
x}\cdot\frac{1}{4\pi}\int_{\mathbb{R}}\mathrm{e}^{-izx}
\hat{\rho}(z)dz. \label{3}
\end{equation}
If $\int_{|x|>1}x^2U(dx)<\infty$ and
$\int_{\mathbb{R}}\left|\psi''(z)+\sigma^2\right|dz<\infty$, then
the Levy measure is absolutely continuous with respect to the
Lebesgue measure, has a bounded continuous density $u(x)$, and
\begin{equation}
u(x)=-\frac{1}{2\pi
x^2}\int_{\mathbb{R}}\mathrm{e}^{-izx}\left(\psi''(z)+\sigma^2\right)dz
\label{4}
\end{equation}}

{\sc{Proof.}} First let us prove formulas (\ref{1}) and (\ref{3}).
Let $f(z)=\log\phi(z)+1/2\cdot\sigma^2z^2$. Therefore, we have
$f(z)=i\gamma z+\int_{\mathbb{R}}
\left(\mathrm{e}^{izx}-1-izxh(x)\right)U(dx)$ and
$f(z)-f(z+\lambda)=\log\frac{\phi(z)}{\phi(z+\lambda)}+1/2\cdot\sigma^2(-\lambda^2-2\lambda
z)=-i\gamma\lambda+\int_{\mathbb{R}}
\left(\mathrm{e}^{izx}-\mathrm{e}^{i(z+\lambda)x}+i\lambda
xh(x)\right)U(dx)$. Let us integrate the both sides of the last
equality over $[-1,1]$ with respect to $\lambda$. On the one hand,
$\hat{\rho}(z)=\int_{-1}^1\left(f(z)-f(z+\lambda)\right)d\lambda=\int_{-1}^1
\log\frac{\phi(z)}{\phi(z+\lambda)}d\lambda-1/3\cdot\sigma^2$. On
the other hand,
$\hat{\rho}(z)=\int_{-1}^1\bigl\{\int_{\mathbb{R}}\bigl(\mathrm{e}^{izx}(1-\mathrm{e}^{i\lambda
x})+ i\lambda xh(x)\bigr)U(dx)\bigr\}d\lambda$. Since
$h(x)\sim1+o(x)$ for $x\to0$ and $h(x)\sim O\left(1/x\right)$ for
$x\to\infty$, it follows that
$\mathrm{e}^{izx}(1-\mathrm{e}^{i\lambda x})+ i\lambda xh(x)\sim
o(x^2)$ for $x\to0$ and $\mathrm{e}^{izx}(1-\mathrm{e}^{i\lambda
x})+ i\lambda xh(x)\sim O(1)$ for $x\to\infty$. Hence we can use
Fubini's theorem and get
$\hat{\rho}(z)=\int_{\mathbb{R}}\left(\mathrm{e}^{izx}\int_{-1}^1(1-\mathrm{e}^{i\lambda
x})d\lambda+ixh(x)\int_{-1}^1\lambda d\lambda\right)U(dx)=
\int_{\mathbb{R}}\mathrm{e}^{izx}2\left(1-\sin x/x\right)U(dx).$
Therefore
$\hat{\rho}(z)=\int_{\mathbb{R}}\mathrm{e}^{izx}\rho(dx)$, where
$\rho(dx)=2\left(1-\sin x/x\right)U(dx)$. This means that the
function $\hat{\rho}(z)$ is the Fourier transform of the measure
$\rho(dx)$. Using the inversion formula for the function
$\hat{\rho}(z)$ [2], we obtain formula (\ref{1}). If
${\int_{\mathbb{R}}\left|\hat{\rho}(z)\right|dz<\infty}$, then the
Levy measure is absolutely continuous with respect to the Lebesgue
measure, has a bounded continuous density $u(x)$, and
$2\left(1-\sin
x/x\right)u(x)=$$1/2\pi\cdot\int_{\mathbb{R}}\mathrm{e}^{-izx}\hat{\rho}(z)dz$.
Thus we have $u(x)=\frac{x}{4\pi(x-\sin
x)}\int_{\mathbb{R}}\mathrm{e}^{-izx}\hat{\rho}(z)dz.$

Now let us prove formulas (\ref{2}) and (\ref{4}). By assumption
$\int_{|x|>1}x^2U(dx)<\infty$, so that, using the monotone
convergence theorem, we get
$\psi''(z)=-\sigma^2+\int_{\mathbb{R}}\left(-x^2\mathrm{e}^{izx}\right)U(dx)$.
It now follows that
$\int_{\mathbb{R}}\mathrm{e}^{izx}x^2U(dx)=-\left(\psi''(z)+\sigma^2\right)$.
Arguing as above, we obtain formulas (\ref{2}) and (\ref{4}).\\

The following examples illustrate the application of the obtained
results for the case of typical infinitely divisible
distributions.

(a) The characteristic function of the {\em normal distribution}
has the form $\phi(z)=\mathrm{e}^{imz-1/2\cdot\sigma^2 z^2}$. It
follows that
$\hat{\rho}(z)=\int_{-1}^1\left(-\sigma^2/2\left(z^2-(z+\lambda)^2\right)-im\lambda\right)
d\lambda-1/3\cdot\sigma^2=0$. Using (\ref{3}), we obtain $u(x)=0$.

(b) The characteristic function of the {\em Cauchy distribution}
has the form $\phi(z)=\mathrm{e}^{-c|z|+i\gamma z}$, $c>0$. It
follows that $\hat{\rho}(z)=
\int_{-1}^1\left(-c|z|+c|z+\lambda|-\gamma\lambda\right)d\lambda=
-2c|z|+c\int_{-1}^1|z+\lambda|d\lambda=\left\{%
\begin{array}{ll}
    c(|z|-1)^2, & \hbox{for $z\in[-1,1]$;} \\
    0, & \hbox{for $z\notin[-1,1]$.} \\
\end{array}%
\right.    $  Using (\ref{3}), we get $u(x)=\frac{cx}{4\pi(x-\sin
x)}\int_{-1}^1\mathrm{e}^{-izx}(|z|-1)^2dz=$$\frac{cx}{4\pi(x-\sin
x)}\frac{4(x-\sin x)}{x^3}=\frac{c}{\pi x^2}$.

(c) The characteristic function of the {\em compound Poisson}
random variable has the form
$\phi(z)=\exp\left\{c\left(\hat{\theta}(z)-1\right)\right\}$,
where $c$ is an intensity of jumps, $\hat{\theta}(z)$ is the
characteristic function of the distribution $\theta(dx)$ of the
jump size. Since the distribution $\theta(dx)$ does not need to
have a density, formula (\ref{1}) is used. We have
$\hat{\rho}(z)=\int_{-1}^1c\left(\hat{\theta}(z)-\hat{\theta}(z+\lambda)\right)d\lambda=$\linebreak$\int_{-1}^1c\left(\int_{\mathbb{R}}\mathrm{e}^{izx}\theta(dx)-
\int_{\mathbb{R}}\mathrm{e}^{i(z+\lambda)x}\theta(dx)\right)d\lambda=$
$2c\int_{\mathbb{R}}\mathrm{e}^{izx}\theta(dx)-
c\int_{\mathbb{R}}\left(\int_{-1}^1\mathrm{e}^{i\lambda
x}d\lambda\right)\mathrm{e}^{izx}\theta(dx)=$\linebreak$
2c\int_{\mathbb{R}}\mathrm{e}^{izx}\left(\frac{x-\sin
x}{x}\right)\theta(dx)$. Using (\ref{1}), we get
$\int_B\left(\frac{x-\sin
x}{x}\right)U(dx)=$\linebreak$\frac{2c}{4\pi}\mathrm{v.p.}\int_{\mathbb{R}}\left(\int_{\mathbb{R}}\mathrm{e}^{izx}
\left(\frac{x-\sin
x}{x}\right)\theta(dx)\int_B\mathrm{e}^{-izy}dy\right)dz=$
$c\int_B\left(\int_{\mathbb{R}} \left(\frac{x-\sin
x}{x}\right)\theta(dx)\frac{1}{2\pi}\mathrm{v.p.}\int_{\mathbb{R}}\mathrm{e}^{iz(x-y)}dz\right)dy=
c\int_B\left(\int_{\mathbb{R}}\left(\frac{x-\sin
x}{x}\right)\theta(dx)\delta(x-y)\right)dy=c\int_B\left(\frac{x-\sin
x}{x}\right)\theta(dx)$. Hence we obtain $\int_B\left(\frac{x-\sin
x}{x}\right)U(dx)=$\linebreak$\int_B\left(\frac{x-\sin
x}{x}\right)\left(c\cdot\theta(dx)\right)$. Since $0<\frac{x-\sin
x}{x}<\infty$ for all $x\neq0$, it now follows that
$U(dx)=c\cdot\theta(dx)$.

(d) The characteristic function of the {\em negative binomial
distribution} has the form \linebreak
$\phi(z)=\frac{p^c}{(1-q\mathrm{e}^{iz})^c},$ where $c>0$,
$p\in(0,1)$, $q=1-p$. We have
$\hat{\rho}(z)=c\int_{-1}^1\log\left(1-q\mathrm{e}^{i(z+\lambda)}\right)d\lambda-
2c\log\left(1-q\mathrm{e}^{iz}\right)$. It is well-known that
$\log(1-x)=-\sum_{k=1}^{\infty}\frac{x^k}{k}$. Let
$x=q\mathrm{e}^{i(z+\lambda)}$, then
$\log(1-q\mathrm{e}^{i(z+\lambda)})=-\sum_{k=1}^{\infty}\frac{q^k\mathrm{e}^{ik(z+\lambda)}}{k}$.
Integrating with respect to $\lambda$ we obtain
$\int_{-1}^1\log(1-q\mathrm{e}^{i(z+\lambda)})d\lambda=-\sum_{k=1}^{\infty}\frac{q^k\mathrm{e}^{ikz}}{k}\cdot
\int_{-1}^1\mathrm{e}^{ik\lambda}d\lambda=
-2\sum_{k=1}^{\infty}\frac{q^k\mathrm{e}^{ikz}}{k}\frac{\sin
k}{k}$. Therefore
$\hat{\rho}(z)=2c\left(-\sum_{k=1}^{\infty}\frac{q^k\mathrm{e}^{ikz}}{k}\frac{\sin
k}{k}+\sum_{k=1}^{\infty}\frac{q^k\mathrm{e}^{ikz}}{k}\right)=2c\cdot\sum_{k=1}^{\infty}
\frac{q^k(k-\sin k)}{k^2}\mathrm{e}^{ikz}$. Using (\ref{1}), we
get $\int_B\left(\frac{x-\sin
x}{x}\right)U(dx)=$\linebreak$\frac{1}{4\pi}\cdot\mathrm{v.p.}\int_{\mathbb{R}}
\left(2c\cdot\sum_{k=1}^{\infty} \frac{q^k(k-\sin
k)}{k^2}\mathrm{e}^{ikz}\right)\left(\int_{B}\mathrm{e}^{-izx}dx\right)dz=$$c\int_B
\left(\sum_{k=1}^{\infty}\frac{q^k(k-\sin
k)}{k^2}\frac{1}{2\pi}\cdot\mathrm{v.p.}\int_{\mathbb{R}}\mathrm{e}^{iz(k-x)}dz\right)dx=
c\sum_{k=1}^{\infty}\frac{q^k(k-\sin k)}{k^2}\int_B\delta(k-x)dx$.
Therefore $\int_B\left(\frac{x-\sin
x}{x}\right)U(dx)=\int_B\sum_{k=1}^{\infty}\left(\frac{k-\sin
k}{k}\right)\frac{cq^k}{k}\delta(k-x)dx$. It now follows that
$U\left(\{k\}\right)=\frac{cq^k}{k}$.

(e) The characteristic function of the {\em hyperbolic cosine
distribution} has the form $\phi(z)=1/\cosh(\pi z/2)$ [1]. Thus we
have $\psi''(z)=\left(\log\phi(z)\right)''=-(\pi^2/4)\phi^2(z)$.
The density of the hyperbolic cosine distribution equals to
$f(x)=1/(\pi\cosh x)$. It is obvious that the convolution product
has the form
$f*f(x)=\int_{\mathbb{R}}f(x-y)f(y)dy=(4/\pi^2)\mathrm{e}^{-x}\int_{\mathbb{R}}
(1+\mathrm{e}^{2y-2x})^{-1}(1+\mathrm{e}^{-2y})^{-1}dy=
\frac{4}{\pi^2}\frac{x}{\mathrm{e}^x-\mathrm{e}^{-x}}$. Therefore
$\int_{\mathbb{R}}\mathrm{e}^{izx}\frac{x}{\mathrm{e}^x-\mathrm{e}^{-x}}dx=
(\pi^2/4)\phi^2(z).$ Using (\ref{4}), we get
$u(x)=\frac{1}{x^2}\cdot\frac{1}{2\pi}
\int_{\mathbb{R}}\mathrm{e}^{-izx}(\pi^2/4)\phi^2(z)dz=
\frac{1}{x^2}\cdot\frac{x}{\mathrm{e}^x-\mathrm{e}^{-x}}=\frac{1}
{x(\mathrm{e}^x-\mathrm{e}^{-x})}$.

(f) The characteristic function of the {\em $\Gamma$-distribution}
has the form $\phi(z)=(1-i\alpha^{-1}z)^{-c}$, $c>0$,
${\alpha>0}$. It follows that
$\psi(z)=\log\phi(z)=-c\log\left(1-i\alpha^{-1}z\right)$,
$\psi''(z)=\frac{c}{\left(i\alpha+z\right)^2}$. Using (\ref{4}),
we obtain
$u(x)=-\frac{1}{x^2}\cdot\frac{1}{2\pi}\int_{\mathbb{R}}\mathrm{e}^{-izx}
\frac{c}{\left(i\alpha+z\right)^2}dz= (c/x)\mathrm{e}^{-\alpha
x}\mathrm{I}_{\{x>\,0\}}$, since
$\int_{\mathbb{R}}\mathrm{e}^{izx}cx\mathrm{e}^{-\alpha
x}\mathrm{I}_{\{x>\,0\}}dx=-\frac{c}{(i\alpha+z)^2}$.\\

The author is grateful to Professor A.N. Shiryaev for constant
attention to this work.

\begin{center}
\bf References
\end{center}
[{\bf1}] W. Feller, {\em An introduction to probability theory and
its applications}, Vol.2, 2nd ed., Wiley, New York, 1971. [{\bf2}]
E. Lukacs, {\em Characteristic functions}, 2nd ed., Griffin,
London, 1970.

$\qquad\qquad\qquad\qquad\qquad\qquad\qquad\qquad\qquad\qquad\qquad\qquad
\qquad\qquad\qquad\qquad\qquad\qquad\qquad\qquad\qquad\qquad\qquad\qquad\qquad$
{\bf E.V. Burnaev}
$\qquad\qquad\qquad\qquad\qquad\qquad\qquad\qquad\qquad\qquad\qquad\;\;\;\;$
{\rm Presented by A.V.
Bulinskiy}\\
Skoltech $\qquad\qquad\qquad\;\;\;\;\;\;\;\;\;\,\,\,\qquad\qquad\qquad\qquad\qquad\qquad\qquad\qquad$ Accepted by editorial board\\
 \\
{\em E-mail}: {\tt e.burnaev@skoltech.ru}

\end{document}